\documentclass[12pt,intlim,righttag]{article}

\usepackage{amscd,amsfonts,amssymb,amsmath,,amsthm,latexsym,array,hhline}

\newcommand{\esp}{\operatornamewithlimits{ess\,sup}}
\newcommand{\la}{\langle}
\newcommand{\ra}{\rangle}
\newcommand{\wt}{\widetilde}

\newcommand{\sg}{\sigma}

\newcommand{\al}{\alpha}
\newcommand{\gm}{\gamma}
\newcommand{\dl}{\delta}
\newcommand{\vp}{\varphi}

\newcommand{\td}{\tilde}
\newcommand{\tht}{\theta}
\newcommand{\ex}{{\cal E}}
\newcommand{\cf}{{\cal F}}
\newcommand{\ift}{\infty}
\newcommand{\be}{\begin}
\newcommand{\ee}{\end}

\font\sm=cmr12 at 6pt

\begin{document}

\title{Solvability of Backward Stochastic Differential Equations with
Quadratic Growth}

\author{Revaz Tevzadze}
\date{~}
\maketitle

\begin{center}

{Georgian--American University, Business School, 3, Alleyway II,
\newline
Chavchavadze Ave. 17\,a,
\newline Georgian Technical Univercity, 77 Kostava str., 0175,
\newline Institute of Cybernetics,  5 Euli str., 0186, Tbilisi,
Georgia
\newline
(e-mail: reztev@yahoo.com) }
\end{center}

\numberwithin{equation}{section}

\begin{abstract}

We prove the existence of the unique solution of a general
Backward Stochastic Differential Equation with quadratic growth
driven by martingales. Some kind of comparison theorem is also
proved.

\bigskip

\noindent {\bf Key words and phrases}:{Backward Stochastic
Differential Equation, Contraction principle, BMO-martingale.}

\noindent
{\bf  Mathematics Subject Classification (2000)}: 90A09, 60H30, 90C39.

\end{abstract}

\

\section{Introduction}

\

In this paper we show a general result of existence and uniqueness
of Backward Stochastic Differential Equation (BSDE) with quadratic
growth driven by continuous martingale. Backward stochastic
differential equations  have been introduced by Bismut \cite{Bs}
for the linear case as equations of the adjoint process in the
stochastic maximum principle. A nonlinear BSDE (with Bellman
generator) was first considered by  Chitashvili \cite{Ch}. He
derived the semimartingale BSDE (or SBE), which can be considered
as a stochastic version of the Bellman equation for a stochastic
control problem, and proved the existence and uniqueness of a
solution. The theory of BSDEs driven by the Brownian motion was
developed by Pardoux and Peng \cite{PP}  for more general
generators. The results of Pardoux and Peng were generalized by
Kobylansky \cite{Kob}, Lepeltier and  San Martin \cite{LS} for
generators with quadratic growth. In the work of Hu at all
\cite{HIm} BMO-martingales were used for BSDE with quadratic
generators in Brownian setting and in \cite{M-T}, \cite{MRT},
\cite{MT7}, \cite{M-S-T22}, \cite{M-T2}, \cite{Mrl} for BSDEs
driven by martingales. By Chitashvili \cite{Ch}, Buckdahn
\cite{B},and El Karoui and Huang \cite{El-H} the well posedness of
BSDE with generators satisfying Lipschitz type conditions was
established. Here we suggest new approach including an existence
and uniqueness of the solution of general BSDE with quadratic
growth. In the earlier papers \cite{M-T}, \cite{MRT}, \cite{MT7},
\cite{M-S-T22}, \cite{M-T2} we studied, as well as Bobrovnytska
and Schweizer \cite{BSc}, the particular cases of BSDE with
quadratic nonlinearities related to the primal and dual problems
of Mathematical Finance. In these works the solutions were
represented as a value function of the corresponding optimization
problems.

The paper is organized as follows. In Section 2 we give some basic
definitions and facts  used in what follows. In Section 3 we show
the solvability of the system of BSDEs for sufficiently small
initial condition and further prove the solvability of one
dimensional BSDE for arbitrary bounded initial data. At the end of
Section 4 we prove the comparison theorem, which generalizes the
results of Mania and Schweizer \cite{m-sc}, and apply this results
to the uniqueness of the solution.

\section{Some basic definitions and assumptions}

\

Let $(\Omega,{\cal F}, {\bf F}=({F}_t)_{t\ge 0},P)$ be filtered
probability space satisfying the usual conditions. We assume that
all local martingales with respect to ${\bf F}$ are continuous.
Here  the time horizon $T<\ift$ is a stopping time and ${\cal F}=F_T$.
Let us
consider Backward Stochastic Differential Equation (BSDE) of the
form
\begin{eqnarray}\label{eq}
dY_t=-f(t,Y_t,\sg_t^*Z_t)dK_t-d\la N\ra_t g_t+Z_t^*dM_t+dN_t,\\
\label{in}
Y_T=\xi
\end{eqnarray}
We suppose that
\begin{itemize}
\item $(M_t,t\ge 0)$ is an $R^n$-valued continuous martingale with
cross-variations matrix $\la M\ra_t=(\la M^i,M^j\ra_t)_{1\le
i,j\le n}$,
\item $(K_t,t\ge 0)$ is a continuous, adapted,
increasing process, such that $\la
M\ra_t=\int_0^t\sg_s\sg^*_sdK_s$ for some predictable, non
degenerate $n\times n$ matrix $\sg$,
\item $\xi$ is ${\cal F}-$measurable an
$R^d$-valued random variable,
\item $f:\Omega\times R^+\times
R^d\times R^{n\times d}\to R^d$ is a stochastic process, such that
for any $(y,z)\in R^d\times R^{n\times d}$ the process
$f(\cdot,\cdot,y,z)$ is predictable,
\item $g:\Omega\times R^+\to
R^{d\times d}$ is a predictable process.
\end{itemize}
The notation $R^{n\times d}$ here denotes the space of ${n\times
d}$-matrix $C$
with Euclidian norm $|C|=\sqrt{{\rm tr}(CC^*)}$. For some stochastic process
$X_t$ and sopping times
$\tau,\;\nu$, such that $\tau\ge\nu$ we denote
$X_{\nu,\tau}=X_\tau-X_\nu$.
For all unexplained notations concerning the martingale theory
used below we refer \cite{J}, \cite{DM} and
\cite{L-Sh}. About {\rm BMO}-martingales  see \cite{DdM} or
\cite{Kaz}.

A solution of the BSDE is a triple $(Y,Z,N)$ of stochastic
processes, such that (\ref{eq}), (\ref{in}) is satisfied and
\begin{itemize}
\item $Y$ is an adapted $R^d$-valued continuous process,
\item $Z$ is an $R^{n\times d}$-valued predictable process,
\item $N$ is an $R^d$-valued continuous martingale, orthogonal to the
basic martingale $M$.
\end{itemize}
One says that $(f,g,\xi)$ is a generator of BSDE
(\ref{eq}),(\ref{in}).

We introduce the following spaces
\begin{itemize}
\item
$L^\ift(R^d)=\{X:\Omega\to R^d, {\cal F}_T-\text{measurable},
||X||_\ift=\underset{\Omega}\esp|X(\omega)|<\infty\}$,
\item
$S^\ift(R^d)=\{\vp:\Omega\times R^+\to R^d,\; \text{continuous,
adapted},||\vp||_\ift=
\underset{[[0,T]]}\esp|\vp(t,\omega)|<\infty\}$,
\item
\be{eqnarray}\label{sp} \notag {H}^2(R^{n\times
d},\sg)=&\{\vp:\Omega\times R^+\to R^{n\times
d},\;\;\text{predictable},\;\;
\\
||\vp||^2_H=\underset{[[0,T]]}\esp
&E(\int_t^T|\sg_s^*\vp_s|^2dK_s|{\cal F}_t)\equiv
\underset{[[0,T]]}\esp E({\rm tr}\la \vp\cdot M\ra_{tT}|{\cal
F}_t\big)<\infty\}, \ee{eqnarray}
\item
${\rm BMO}(Q)=\{ N, \; R^d-\text{valued}\; Q-\text{martingale}\;
||N||_Q^2=\underset{[[0,T]]}\esp E^Q(tr\la N\ra_{tT}|{\cal
F}_t)<\infty\}$
\end{itemize}
We also use  the notation $|r|_{2,\ift}$
for the norm  $||\int_0^Tr_s^2dK_s||_\ift$.

\noindent
The norm of the triple is defined as
$$
||(Y,Z,N)||^2=||Y||^2+||Z||_H^2+||N||_P^2.
$$
Throughout the paper we use the condition

A) There exist a constant $\theta$ and predictable processes
$$
\alpha:\Omega\times R^+\to R^d,\;\Gamma:\Omega\times R^+\to
Lin(R^{n\times d},R^d),\;\; r:\Omega\times R^+\to R,
$$
such that the following conditions $\int_0^T r_sdK_s,\;\int_0^T r_s^2dK_s\in
L^\ift,\;\Gamma(\sg^{-1})\in H_T^2,$  $|\alpha_t|\le r_t,|g_t|\le
\theta^2$ and \be{eqnarray}\label{est}
|f(t,y_1,z_1)-f(t,y_2,z_2)-\alpha_t(y_1-y_2)-\Gamma_t(z_1-z_2)|\\
\notag \le(r_t|y_1-y_2|+
\theta|z_1-z_2|)(r_t(|y_1|+|y_2|)+\theta(|z_1|+|z_2|)).
\ee{eqnarray} are satisfied.

Sometimes we use the more restrictive conditions \be{itemize}
\item[B1)] $\int_0^T|f(t,0,0)|dK_t+|g_t|\le \theta^2$ for all
$t\in[0,T]$,

\item[B2)] $|f_y(t,y,z)|\le r_t,\; |f_z(t,y,z)|\le r_t+\tht|z|$ for all
$(t,y,z)$,

\item[B3)] $|f_{yy}(t,y,z)|\le r_t^2,\;|f_{yz}(t,y,z)|\le \tht
r_t,\;|f_{zz}(t,y,z)|\le \tht^2$
for all $(t,y,z)$. \ee{itemize}

{\bf Remark 1.} Condition A) follow from conditions B1)-B3), since
using notations $\delta y=y_1-y_2,\;\delta z=z_1-z_2$ for
$\alpha_t=f_y(t,0,0), \; \Gamma_t= f_z(t,0,0)$ by the mean value
theorem we have
$$
|f(t,y_1,z_1)-f(t,y_2,z_2)-\al_t\dl y-\Gamma_t(\dl z)|
$$
$$
=|f_y(t,\nu y_1+(1-\nu)y_2,\nu z_1+(1-\nu)z_2)\dl y-f_y(t,0,0)\dl y|
$$
$$
+f_z(t,\nu y_1+(1-\nu)y_2,\nu z_1+(1-\nu)z_2)(\dl z)
-f_z(t,0,0)(\dl z)|,
$$
for some $\nu\in [0,1].$ Using again mean value theorem we obtain
that
$$
|f(t,y_1,z_1)-f(t,y_2,z_2)-\al_t\dl y-\Gamma_t(\dl z)|
$$
$$
\le (|\nu y_1+(1-\nu)y_2|\max_{y,z}|f_{yy}(t,y,z)|+|\nu z_1+(1-\nu)z_2
|\max_{y,z}|f_{yz}(t,y,z)|)|\dl y|
$$
$$
+(|\nu y_1+(1-\nu)y_2|\max_{y,z}|f_{yz}(t,y,z)|+|\nu
z_1+(1-\nu)z_2|\max_{y,z}|f_{zz}(t,y,z)|)|\dl z|
$$
$$
\le [r_t^2(|y_1|+|y_2|)+r_t\tht(|z_1|+|z_2|)]|\dl y|
+[r_t\tht(|y_1|+|y_2|)+\tht^2(|z_1|+|z_2|)]|\dl z|
$$
$$
=(r_t|\dl y|+\tht|\dl z|)(r_t(|y_1|+|y_2|)+\tht(|z_1|+|z_2|).
$$

{\bf Remark 2.} If $d=1$ the operator $\Gamma_t$ is given by an
$n-$dimensional vector $\gm_t$ such that $\Gamma_t(z)=\gm_t^*z$.
Thus inequality in A) can be rewritten as
$$
|f(t,y_1,z_1)-f(t,y_2,z_2)-\alpha_t\dl y-\gm_t^*\dl z|
$$
$$
\le(r_t|\dl y|+
\theta|\dl z|)(r_t(|y_1|+|y_2|)+\theta(|z_1|+|z_2|)).
$$

The main statement of the paper is the following

{\bf Theorem 1.} Let $\xi\in L^\infty,\;d=1$ and conditions B1)-B3) are
satisfied. Then there exists
a unique triple $(Y,Z,N)$, where
$Y\in S^\ift,Z\in H^2, N\in BMO$, that satisfies
equation (\ref{eq}),(\ref{in}).

\section{Existence of the solution}

\

First we prove the existence and uniqueness of the solution for a
sufficiently small initial data.

{\bf Proposition 1}. Let $f$ and $g$ satisfy condition $A)$ with
$\al=0$ and $\gm_t=0$. Then for $\xi$ with the norm
$||\xi||_\ift<\frac{1}{32\beta},
\;\beta=8\max(|r|_{2,\ift}^2,\tht^2)$ there exists a unique
solution $(Y,Z,N)$ of BSDE \be{eqnarray}\label{eqq}
dY_t=(f(t,0,0)-f(t,Y_t,\sg_t^*Z_t))dK_t+d\la N\ra_t
g_t+Z_t^*dM_t+dN_t,\\
\notag
Y_T=\xi,
\ee{eqnarray}
with the norm $||(Y,Z,N)||\le R$, where $R$ is a constant satisfying
the inequality
$4||\xi||^2_\ift+\beta^2R^4\le R^2$, namely $R=2\sqrt2||\xi||_\ift.$

Moreover if $||\xi||_\ift+||\int_0^\ift|f(s,0,0)|dK_s||_\ift$ is
small enough then BSDE (\ref{eq}) admits a unique solution.

{\it Proof}. We define the mapping $(Y,Z,N)=F(y,z,n),\;\; n\;\;
\text{is orthogonal to}\;\;\; M,\;\;\;\\(y,z\cdot M+n)\in
S_T^\ift\times BMO(P)$ by the relation \be{eqnarray}\label{eq1}
\notag dY_t=(f(t,0,0)-f(t,y_t,\sg_t^*z_t))dK_t+d\la
n\ra_tg_t+Z_t^*dM_t+dN_t,\\
Y_T=\xi.
\ee{eqnarray}
Using the Ito formula for $|Y_t|^2$ we obtain that
\be{eqnarray}\label{ito}
\notag
|Y_t|^2=|\xi|^2+2\int_t^TY_s^*(f(s,y_s,\sg_s^*z_s)-f(s,0,0))dK_t\\
\notag
+2\int_t^TY_s^*d\la n\ra_sg_s-\int_t^T{\rm tr}Z_s^*d\la M\ra_sZ_s
-{\rm tr}\la N\ra_{tT}-\int_t^TY_s^*Z_s^*dM_s-\int_t^TY_s^*dN_{s}.
\ee{eqnarray}
If we take the conditional expectation and use (\ref{sp}) and the
elementary inequality $2ab\le\frac{1}{4}a^2+4b^2$ we get
\be{eqnarray}\label{itt}
\notag
|Y_t|^2+E(\int_t^T|\sg_s^*Z_s|^2dK_s+tr\la N\ra_{tT}|{\cal F}_t)\le
||\xi||^2+\frac{1}{4}||Y||_\ift^2\\
+4E^2(\int_t^T|f(s,y_s,\sg_s^*z_s)-f(s,0,0)|dK_s+
\int_t^T|g_s|d{\rm tr}\la n\ra_s|{\cal F}_t). \ee{eqnarray} Thus
using condition A), identities \be{eqnarray}\label{idd} {\rm
tr}\la z\cdot M\ra_t={\rm tr}\int_0^tz_s^*d\la M\ra_sz_s=
\int_0^t{\rm
tr}(z_s^*\sg_s\sg_s^*z_s)dK_s=\int_0^t|\sg_s^*z_s|^2dK_s
\ee{eqnarray} and explicit inequalities \be{eqnarray}\label{12}
\notag \frac{1}{2}(||Y||_\ift^2+||Z\cdot M+N||_{\text{\sm BMO}}^2)
\le \max(||Y||_\ift^2,||Z\cdot M+N||_{\text{\sm BMO}}^2)\\
\notag \le
\underset{[[0,T]]}\esp[|Y_t|^2+E(\int_t^T|\sg_s^*Z_s|^2dK_s+tr\la
N\ra_{tT}|{\cal F}_t)] \ee{eqnarray}
we obtain from (\ref{itt})
\be{eqnarray}\label{12} \notag
\frac{1}{4}||Y||_\ift^2+\frac{1}{2}||Z\cdot M+N||_{\text{\sm BMO}}^2\le
||\xi||^2\\
\notag
+4\underset{[[0,T]]}\esp
E^2(\int_t^T|f(s,y_s,\sg_s^*z_s)-f(s,0,0)|dK_t+
\tht^2{\rm tr}\la n\ra_{tT}|{\cal F}_t)\\
\notag
\le ||\xi||^2+16\underset{[[0,T]]}\esp
E^2(\int_t^Tr_s^2y_s^2dK_s+\tht^2{\rm tr}\la z\cdot M+n\ra_{tT}|{\cal F}_t)\\
\notag
\le ||\xi||^2+16|r|_{2,\ift}^4||y||_\ift^4+16\tht^4||z\cdot
M+n||^4_{\text{\sm BMO}}.
\ee{eqnarray}
Therefore
\be{eqnarray*}\label{13} \be{gathered}
||Y||_\ift^2+||Z\cdot M+N||_{\text{\sm BMO}}^2\le 4||\xi||^2\\
+64|r|_{2,\ift}^4||y||^4_\ift+64\tht^4||z\cdot M+n||_{\text{\sm BMO}}^4\\
\le 4||\xi||^2+\beta^2(||y||_\ift^2+||z\cdot M+n||^2_{\text{\sm
BMO}})^2,
\ee{gathered}
\ee{eqnarray*}
where
$\beta=8\max(|r|_{2,\ift}^2,\tht^2)$. We can pick $R$ such that
$$
4||\xi||^2+\beta^2R^4\le R^2
$$
if and only if
$||\xi||_\ift \le\frac{1}{4\beta}$. For instance
$R=2\sqrt2||\xi||_\ift$ satisfies this quadratic inequality. Therefore
the ball
$$
{\cal B}_R=\{(Y,Z\cdot M+N)\in S^\ift\times{\rm BMO},\;N\bot
M,\;\; ||Y||_\ift^2+||Z\cdot M+N||_{\text{\sm BMO}}^2\le R^2\}
$$
is such that
$F({\cal B}_R)\subset{\cal B}_R$.

Similarly for $(y^j,z^j\cdot M+n^j)\in{\cal B}_R,\;j=1,2$ using the
notations
$\delta y=y^1-y^2,\;\delta z=z^1-z^2, \;\dl n=n^1-n^2$ we can show
that
\be{eqnarray*}\label{14}
\be{gathered}
||\dl Y||_\ift^2+||\dl Z\cdot M+\dl N||_{\text{\sm BMO}}^2\\
\le 4\underset{[[0,T]]}\esp
E^2\big(\int_t^T|f(s,y_s^1,\sg_s^*z_s^1)-
f(s,y^2_s,\sg_s^*z_s^2)|dK_s\\
+\int_t^T|g_s|d{\rm var}({\rm tr}\la\dl n,n^1+n^2\ra)_s|{\cal F}_t)\\
\le 8\underset{[[0,T]]}\esp E\big(\int_t^T(r_s^2|\dl y_s|^2+
\tht^2|\sg_s^*\dl z_s|^2dK_s)|{\cal F}_t)\\
\times
E\big(\int_t^T(r_s(|y_s^1|+|y_s^2|)
+\tht(|\sg_s^*z^1_s|+|\sg_s^*z_s^2|))^2dK_s)|{\cal F}_t)\\
+\tht^2E({\rm tr}\la\dl n\ra_{tT}|\cf_s)
E({\rm tr}\la n^1+n^2\ra_{tT}|\cf_t)
\ee{gathered}
\ee{eqnarray*}
Again using the equalities ({\ref{idd})
we can pass to the norm. Thus
\be{eqnarray*}\label{114}
\be{gathered}
||\dl Y||_\ift^2+||\dl Z\cdot M+\dl N||_{\text{\sm BMO}}^2\\
\le 8(|r|_{2,\ift}^2||\dl y||^2_\ift+\tht^2||\dl z\cdot M||_{\text{\sm
BMO}}^2)\\
\times(|r|_{2,\ift}^2(||y^1||^2_\ift+||y^2||^2_\ift)
+\tht^2(||z^1\cdot M||_P^2+||z^2\cdot M||_P^2)\\
+2\tht^2||\dl n||^2_{\text{\sm BMO}}
(||n^1||^2_{\text{\sm BMO}}+||n^2||^2_{\text{\sm BMO}})^2).
\ee{gathered}
\ee{eqnarray*}
Since $||z^1\cdot M||,||z^2\cdot M||\le R,||n^1||,||n^2||\le R$ we get
\be{eqnarray}\label{144}
\be{gathered}
||\dl Y||_\ift^2+||\dl Z\cdot M+\dl N||_{\text{\sm BMO}}^2\\
\le128\beta^2 R^2(||\dl y||_\ift^2+||\dl z\cdot M||^2_{\text{\sm
BMO}})
+4\beta^2R^2||\dl n||^2_{\text{\sm BMO}}\\
\le128\beta^2 R^2(||\dl y||_\ift^2+||\dl z\cdot M+\dl n||^2_{\text{\sm
BMO}}).
\ee{gathered}
\ee{eqnarray}

Now we can take $R=2\sqrt2||\xi||_\ift< \frac{1}{8\sqrt 2\beta}$. This
means that $||\xi||_\ift< \frac{1}{32\beta}$ and
$F$ is contraction on ${\cal B}_R$. By contraction
principle the mapping $F$ admits a unique fixed point, which is
the solution of (\ref{eqq}). \qed

From now we suppose that $d=1$.

{\bf Lemma 1}. Let condition A) is satisfied. Then the
generator $(\bar f,\bar g,\bar\xi)$, where
$$
\bar f(t,\bar y,\bar
z)=e^{\int_0^t\al_sdK_s}(f(t,e^{-\int_0^te^{\al_sdK_s}}
\bar y,e^{-\int_0^te^{\al_sdK_s}}\bar z)-f(t,0,0))-\al_t\bar
y-\gm_t^*\bar z,
$$
$$
\bar g_t=e^{-\int_0^t\al_sdK_s}g_t\;\;\; {\rm and}\;\;\; \bar
\xi=e^{\int_0^T\al_sdK_s}\xi,\;\;
$$
satisfies  condition A) with $\al=0,\;\gm=0,\;\bar
r_t=r_te^{||\int_0^\ift r_sdK_s||_\ift},\; \text{and}\;\bar\tht=\tht
e^{||\int_0^T r_sdK_s||_\ift}$.

Moreover, $(Y,Z,N)$ is a solution of BSDE (\ref{eqq}) if and only
if
$$
(\bar Y_t,\bar Z_t,\bar N_t)=
(e^{\int_0^t\al_sdK_s}Y_t,e^{\int_0^t\al_sdK_s}Z_t,\int_0^te^{\int_0^s\al_udK_u}dN_s)
$$
is a solution w.r.t. measure $d\bar P=\ex_T((\gm\sg^{-1})\cdot
M)dP$ of BSDE
\begin{eqnarray}\label{beq}
d\bar Y_t=-\bar f(t,\bar Y_t,\sg_t^*\bar Z_t)dK_t-d\la \bar N\ra_t \bar
g_t+
\bar Z_t^*d\bar M_t+d\bar N_t,\\
\notag \bar Y_T=\bar \xi,
\end{eqnarray}
 where
$\bar M_t=M_t-\la(\gm\sg^{-1})\cdot M,M\ra_t$.

{\it Proof}.
Condition A) for $(\bar f,\bar g,\bar\xi)$ is satisfied since by
(\ref{est})
$$
|\bar f(t,\bar y_1,\bar z_1)-\bar f(t,\bar y_2,\bar z_2)|
$$
$$
\le e^{\int_0^t\al_sdK_s}(r_t|\dl \bar y|+ \theta|\dl \bar
z|)(r_t(|\bar y_1|+|\bar y_2|)+\theta(|\bar z_1|+|\bar z_2|))
$$
$$
\le(\bar r_t|\dl \bar y|+
\bar \theta|\dl \bar z|)(\bar r_t(|\bar y_1|+|\bar
y_2|)+\bar\theta(|\bar z_1|+|\bar z_2|)).
$$
On the other hand
using the Ito formula we have
\be{eqnarray*}\label{15}
\be{gathered}
d\bar
Y_t=e^{\int_0^t\al_sdK_s}dY_t+\al_te^{\int_0^t\al_sdK_s}Y_tdK_t\\
=e^{\int_0^t\al_sdK_s}(f(t,0,0)-f(t,Y_t,\sg_t^*Z_t))dK_t
+e^{\int_0^t\al_sdK_s}d\la N\ra_tg_t\\
+e^{\int_0^t\al_sdK_s}Z_t^*dM_t+e^{\int_0^t\al_sdK_s}dN_t+\al_t\bar
Y_tdK_t
\ee{gathered}
\ee{eqnarray*}
Taking into account that
\be{eqnarray*}\label{155}
\be{gathered}
e^{\int_0^t\al_sdK_s}(f(t,0,0)-f(t,Y_t,\sg_t^*Z_t))+\al_t\bar Y_t\\
=-\bar f(t,\bar Y_t,\sg_t^*\bar Z_t)-\gm_t\sg_t^*\bar Z_t,\\
e^{\int_0^t\al_sdK_s}d\la N\ra_tg_t=d\la\bar
N\ra_te^{-\int_0^t\al_sdK_s}g_t=d\la\bar N\ra_t\bar g_t
\ee{gathered} \ee{eqnarray*} and
$$
\bar Z\cdot M-\int_0^\cdot\gm_t\sg_t^*\bar Z_tdK_t=\bar Z\cdot M-
\int_0^\cdot\gm_t\sg_t^{-1}\sg_t\sg_t^*\bar Z_tdK_t
$$
$$
= \bar Z\cdot M-\int_0^\cdot\gm_t\sg_t^{-1}d\la M\ra_t\bar Z_t =
\bar Z\cdot M-\la(\gm\cdot\sg^{-1})\cdot M,\bar Z\cdot M\ra =\bar
Z\cdot \bar M
$$
we obtain
$$
d\bar Y_t=-\bar f(t,\bar Y_t,\sg_t^*\bar Z_t)dK_t-d\la\bar N\ra_t\bar
g_t+\bar Z_td\bar M_t+
d\bar N_t.
$$
Here $\bar M$ is a local martingale w.r.t. $\bar P$ by Girsanov
theorem.

{\bf Corollary 1}. Let $f$ and $g$ satisfy condition A) and
$||\xi||_\ift\le\frac{1}{32\beta}\exp(-2||\int_0^T
r_sdK_s||_\ift)$. Then there exist the solution of (\ref{eqq})
with the norm $||Y||_\ift^2+||Z\cdot\bar M+N||_{\rm BMO(\bar
P)}^2\leq \frac{1}{128\beta^2}.$

{\it Proof}. Obviously that $$||Y||_\ift^2+||Z\cdot\bar M+N||_{\rm BMO(\bar
P)}^2\leq
\left(||\bar Y||_\ift^2+||\bar Z\cdot\bar M+\bar N||_{\rm BMO(\bar
P)}^2\right)\exp(2||\int_0^T r_sdK_s||_\ift)$$
$$\leq 8||\bar\xi||_\ift^2\exp(2||\int_0^T r_sdK_s||_\ift)
\leq 8||\xi||_\ift^2\exp(4||\int_0^T r_sdK_s||_\ift).$$
From $||\xi||_\ift\le\frac{1}{32\beta}\exp(-2||\int_0^T r_sdK_s||_\ift)$ follows that
$8||\xi||_\ift\exp(4||\int_0^T r_sdK_s||_\ift)\leq \frac{1}{128\beta^2}.$
Hence we get $||Y||_\ift^2+||Z\cdot\bar M+N||_{\rm BMO(\bar P)}^2\leq\frac{1}{128\beta^2}.$

{\bf Corollary 2.} Let generator $(f,g,\xi)$ satisfies conditions
B1)-B3) and $(\tilde Y_t,\tilde  Z_t,\tilde  N_t)$ be a solution
of (\ref{eqq}). Then BSDE
\begin{eqnarray}\label{heq}
d\hat Y_t=(f(t,\tilde Y_t,\sg_t^*\tilde Z_t)-
f(t,\hat Y_t+\tilde Y_t,\sg_t^*\hat Z_t+\sg_t^*\tilde Z_t))dK_t\\
\notag -d(\la \hat N\ra_t+2\la \tilde N,\hat N\ra_t)g_t+
\hat Z_t^*dM_t+d\hat N_t,\\
\notag
\hat Y_T=\hat \xi
\end{eqnarray}
satisfy condition A) with $-\hat f(t,y,z)=f(t,\tilde
Y_t,\sg_t^*\tilde Z_t)- f(t,y+\tilde Y_t,z+\sg_t^*\tilde Z_t)$,
$\al_t=f_y(t,\tilde Y_t,\sg_t^*\tilde Z_t)$, $\gm_t=f_z(t,\tilde
Y_t,\sg_t^*\tilde Z_t)$ and the new probability measure ${\mathcal
E}_T(2g\cdot \tilde N)dP$. Moreover (\ref{heq}) admits a
unique solution $(\hat Y_t,\hat Z_t,\hat N_t)$ if
$||\hat\xi||_\ift\le\frac{1}{32\beta}\exp(-2||\int_0^\cdot
r_sdK_s||_\ift)$.

{\it Proof}. Using a change of measure the equation (\ref{heq}) reduces to
equation of type (\ref{eqq}). By previous corollary we obtain the existence
and uniqueness of the BSDE.

{\bf Lemma 2}. Let conditions B1)-B3) be satisfied and
random variables $\tilde \xi$ and $\hat\xi$ be such that
$max(||\tilde\xi||_\ift,||\hat\xi||_\ift)\le\frac{1}{32\beta}e^{-2||\int_0^T
r_s^2dK_s||_\ift}.$ Then there exist solutions of  BSDEs
(\ref{heq}) and
\begin{eqnarray}\label{16}
d\tilde Y_t=(f(t,0,0)-
f(t,\tilde Y_t,\sg_t^*\tilde Z_t))dK_t
-d\la \tilde N\ra_tg_t+
\tilde Z_t^*dM_t+d\tilde N_t,\\
\notag
\tilde Y_T=\tilde \xi
\end{eqnarray}
 and the triple $(Y,Z,N)=(\tilde Y+\hat Y,\tilde
Z+\hat Z,\tilde N+\hat N)$ satisfies  BSDE
\begin{eqnarray}\label{18}
\notag dY_t=(f(t,0,0)- f(t,Y_t,\sg_t^* Z_t))dK_t -d\la N\ra_tg_t+
Z_t^*dM_t+dN_t,\\
\notag
Y_T=\tilde \xi+\hat\xi.
\end{eqnarray}

{\it Proof}. Similarly to the Remark from Section 1 we can show
that for $\hat f(t,y,z)= f(t,\tilde Y_t,\sg_t^*\tilde Z_t)-
f(t,y+\tilde Y_t,\sg_t^*z+\sg_t^*\tilde Z_t), \;
\al_t=f_y(t,\tilde Y_t,\sg_t^*\tilde Z_t),\; \gm_t=f_z(t,\hat
Y_t,\sg_t^*\hat Z_t)$  the estimate
$$
|\hat f(t,y_1,z_1)-\hat f(t,y_2,z_2)-\al_t\dl y-\gm_t^*\dl z|
$$
$$
\le (r_t|\dl y|+\tht|\dl z|)(r_t(|y_1|+|y_2|)+\tht(|z_1|+|z_2|)).
$$
holds.

Now  by Lemma 1 and Corollary 2 of Lemma 1 we obtain the
solvability of both equations (\ref{16}),(\ref{heq}). \qed

{\bf Proposition 2.} Let $f$ and $g$ satisfy condition B1)-B3) and
$\xi\in L^\ift$. Then BSDE (\ref{eq}) admits a solution
$(Y,Z\cdot M+N)\in S^\ift\times{\rm BMO}$.

{\it Proof}. An arbitrary $\xi\in L^\ift(R)$ can be represented as sum
$\xi=\sum_{i=1}^m\xi_i$ with
$||\xi_i||_\ift\le\frac{1}{32\beta}\exp(-2||\int_0^\cdot
r_sdK_s||_\ift)$.
Denote by $(Y^j,Z^j,N^j),\;j=1,...,m$ the solution of
\begin{eqnarray}\label{ieq}
\notag
dY^j_t=(f(t,Y^0_t+...+Y^{j-1}_t,\sg_t^*(Z^0_t+...+Z^{j-1}_t))\\
-f(t,Y^0_t+...+Y^j_t,\sg_t^*(Z^0_t+...+Z^j_t))dK_t\\
\notag -d(\la N^j\ra_t+2\la N^j,N^0+...+N^{j-1}\ra_t)g_t+
Z_t^{j*}dM_t+dN^j_t,\\
\notag Y^j_T=\xi^j\\
\notag Y^0=0,\;\;Z^0=0\;\;N^0=0.
\end{eqnarray}
By Corollary 1 we get
$$
||Y^j||_\ift^2+||Z^j\cdot M^j+N^j||_{{\rm BMO(P^j)}}^2\leq
\frac{1}{128\beta^2},
$$
where $dP^j={\cal E}_T(\int_0^\cdot
f_z(s,Y^0_s+...Y_s^{j-1},\sg_s^*(Z_s^0+...+Z_s^{j-1}))\sg_s^{-1}dM_s)dP,
$ and $ M^j=M-\la
f_z(\cdot,Y^0+...+Y^{j-1},\sg^*(Z^0+...+Z^{j-1}))\sg^{-1}\cdot
M,M\ra.$

Using Lemma 2 we get the existence of a solution for BSDE
\begin{eqnarray*}
d\bar Y_t=(f(t,0,0)-f(t,\bar Y_t,\sg_t^*Z_t))dK_t-d\la N\ra_tg_t+
Z_t^*dM_t+dN_t,\\
\notag \bar Y_T=\xi.
\end{eqnarray*}
Since $\int_0^Tf(t,0,0)dK_t$ is bounded we can apply the above
argument with $f$ replaced by $\bar
f(t,y,z)=f(t,y-\int_0^tf(s,0,0)dK_s,z)$ to get the existence of
solution
\begin{eqnarray*}
d\bar Y_t=(f(t,0,0)-f(t,\bar
Y_t-\int_0^tf(s,0,0)dK_s,\sg_t^*Z_t))dK_t-d\la N\ra_tg_t+
Z_t^*dM_t+dN_t,\\
\notag \bar Y_T=\xi+\int_0^Tf(s,0,0)dK_s.
\end{eqnarray*}
Obviously $Y_t=\bar Y_t-\int_0^tf(s,0,0)dK_s$ is a solution of
BSDE (\ref{eq}),(\ref{in}).

\

\section{A comparison theorem for BSDEs}

\

Let us consider BSDE (\ref{eq}),(\ref{in}) in the case $d=1$.

{\bf Lemma 3.} Let $\xi\in L_\infty$ and assume that
there are positive constants $C(f),C(g)$, increasing function $\lambda:
R^+\to R^+$,
bounded on all bounded subsets and a predictable
process $k\in H^2(R,1)$ such that
\begin{equation}\label{2.3}
|f(t,y,z)|\le k^2_t\lambda(|y|)+ C(f)z^2,
\end{equation}
\begin{equation}\label{2.4}
|g(t)|\le C(g).
\end{equation}
Then the martingale part of any bounded solution of
(\ref{eq}),(\ref{in}) belongs to the space ${\rm BMO}(P)$.

{\it Proof.} Let $Y$ be a solution of (\ref{eq}),(\ref{in}) and there
is a
constant $C>0$ such that
$$
|Y_t|\le C\;\;\;\text{a.s}\;\;\;\text{for all}\;\;\; t.
$$
Applying the It\^o formula for
$\exp\{\beta Y_T\}-\exp\{\beta Y_\tau\}$ and using the
boundary condition $Y_T=\xi$ we have
\be{eqnarray}\label{3.5}
\be{gathered}
\frac{\beta^2}{2}\int_\tau^Te^{\beta Y_s}Z^*_sd\la M\ra_sZ_s+
\frac{\beta^2}{2}\int_\tau^Te^{\beta Y_s}d\la N\ra_s\\
-\beta\int_\tau^Te^{\beta Y_s}f(s,Y_s,Z_s)dK_s-
\beta\int_\tau^Te^{\beta Y_s}g(s)d\la N\ra_s\\
+\beta\int_\tau^Te^{\beta Y_s}Z_s^*dM_s+
\beta\int_\tau^Te^{\beta Y_s}dN_s=
e^{\beta \xi}-e^{\beta Y_\tau}\le e^{\beta C},
\ee{gathered}
\ee{eqnarray}
where $\beta$ is a constant yet to be determined.

If $Z\cdot M$ and $N$ are square integrable martingales taking
conditional expectations in (\ref{3.5}) we obtain
\be{eqnarray*}\label{35} \be{gathered}
\frac{\beta^2}{2}E\big(\int_\tau^Te^{\beta Y_s}Z^*_sd\la
M\ra_sZ_s|F_\tau\big) +\frac{\beta^2}{2}E\big(\int_\tau^Te^{\beta
Y_s}d\la N\ra_s|F_\tau\big)
\\
\le e^{\beta C}+\beta E\big(\int_\tau^Te^{\beta Y_s}
|f(s,Y_s,Z_s)|dK_s|F_\tau\big)
+\beta E\big(\int_\tau^Te^{\beta Y_s}|g(s)|d\la N\ra_s|F_\tau\big)
\ee{gathered}
\ee{eqnarray*}
Now if we use the estimates (\ref{2.3}),(\ref{2.4}) we get
\be{eqnarray*}\label{335}
\be{gathered}
\frac{\beta^2}{2}E\big(\int_\tau^Te^{\beta Y_s}Z^*_sd\la
M\ra_sZ_s|F_\tau\big)
+\frac{\beta^2}{2}E\big(\int_\tau^Te^{\beta Y_s}d\la
N\ra_s|F_\tau\big)\\
\le e^{\beta C}+\beta\lambda(C)E\big(\int_\tau^Te^{\beta
Y_s}k_s^2dK_s|F_\tau\big)\\
+\beta C(f)E\big(\int_\tau^Te^{\beta
Y_s}|\sg_s^*Z_s|^2dK_s|F_\tau\big)
+\beta E\big(\int_\tau^Te^{\beta Y_s}|g(s)|d\la N\ra_s|F_\tau\big)\\
\le e^{\beta C}+\beta\lambda(C)E\big(\int_\tau^Te^{\beta
Y_s}k_s^2dK_s|F_\tau\big)\\
+\beta C(f)E\big(\int_\tau^Te^{\beta Y_s}|Z_s^*d\la
M\ra_sZ_s|^2|F_\tau\big)
+C(g)\beta E\big(\int_\tau^Te^{\beta Y_s}d\la N\ra_s|F_\tau\big).
\ee{gathered}
\ee{eqnarray*}
Conditions (\ref{2.3}) and (\ref{2.4}) imply that
\be{eqnarray}\label{3.6}
\be{gathered}
(\frac{\beta^2}{2}-\beta C(f))E\big(\int_\tau^T
e^{\beta Y_s}Z^*_sd\la M\ra_sZ_s|F_\tau\big)+
\\
+(\frac{\beta^2}{2}-\beta C(g))E\big(\int_\tau^T
e^{\beta Y_s}d\la N\ra_s|F_\tau\big)\le
\\
\le e^{\beta C}+\beta\lambda(C) E\big(\int_\tau^Te^{\beta Y_s}
k^2_sdK_s|F_\tau\big).
\ee{gathered}
\ee{eqnarray}
Taking $\beta = 4 \overline C$, where $\overline C=max(C(f),C(g))$,
from (\ref{3.6}) we have
$$
4\overline C^2[E\big(\int_\tau^T
e^{\beta Y_s}Z^*_sd\la M\ra_sZ_s|F_\tau\big)+
E\big(\int_\tau^T
e^{\beta Y_s}d\la N\ra_s|F_\tau\big)]\le
$$
$$
\le e^{4 C \overline C}\big(4\overline C\lambda(C)||k||_{H}+1\big).
$$
Since $Y\ge -C$, from the latter inequality we finally obtain the
estimate
$$
E\big(\la Z\cdot M\ra_{\tau T}|F_\tau\big)+
E\big(\la N\ra_{\tau T}|F_\tau\big)\le
$$
\begin{equation}\label{3.7}
\le \frac{e^{8 C\overline C}[4\overline C\lambda(C)||k||_{H}+1]}
{4 \overline C^2}
\end{equation}
for any stopping time $\tau$, hence
$Z\cdot M, N\in BMO$.

For general $Z\cdot M$ and $N$ we stop at $\tau_n$ and derive
(\ref{3.7}) with $T$ replaced $\tau_n$. Letting $n\to\infty$ then
completes the proof. \qed

Further we use some notations. Let $(Y,Z),(\wt Y,\wt Z)$ be two pairs
of processes and $(f,g,\xi),(\td f,\td g,\td\xi)$ two triples of
generators. Then we denote:
$$
\dl f=f-\td f,\;\;\dl g=g-\td g,\;\;\dl\xi=\xi-\td\xi,
$$
$$
\partial_yf(t,Y_t,\wt Y_t,Z_t)\equiv\partial f_y(t)
=\frac{f(t,Y_t,Z_t)-f(t,\wt Y_t,Z_t)}{Y_t-\wt Y_t}
$$
$$
\text{for all}\;\;j=1,..,n,\;\;\partial_{j}f(t,\wt Y_t,Z_t,\wt
Z_t)\equiv\partial_j f(t)
$$
$$
=\frac{
f(t,\wt Y_t,Z^1_t,...,Z^{j-1}_t,Z^j_t,\wt Z^{j+1}_t,...,\wt Z^n_t)-
f(t,\wt Y_t,Z^1_t,...,Z^{j-1}_t,\wt Z^j_t,\wt Z^{j+1}_t,...,\wt
Z^n_t)}
{Z^j_t-\wt Z^j_t},
$$
$$
\nabla f(t)=(\partial_{1}f(t),...,\partial_{n}f(t))^*
$$
Thus we have
\be{eqnarray}\label{38}
f(t,Y_t,Z_t)-f(t,\wt Y_t,\wt Z_t)=\partial_yf(t)\dl Y_t+\nabla
f(t)^*\dl Z_t.
\ee{eqnarray}

{\bf Theorem 2.} Let $Y$ and $\wt Y$ be the bounded solutions of
SBE (\ref{eq}) with generators $(f,g,\xi)$ and $(\td f,\td
g,\td\xi)$ respectively, satisfying the conditions of Lemma 3.

If $\xi\ge \td\xi$ (a.s), $f(t,y,z)\ge \td f(t,y,z)$
($\mu^{K}$-a.e.), $g(t)\ge \td g(t)$ ($\mu^{\la N\ra}$-a.e.) and
$f$ (or $\td f$) satisfies the following Lipschitz condition:

L1) for any $Y,\wt Y,Z$
$$
\frac{f(t,Y_t,Z_t)-f(t,\wt Y_t,Z_t)}{Y_t-\wt Y_t}\in S^\ift,
$$

L2) for any $Z,\wt Z\in H^2$ and any bounded process $Y$
$$
(\sg_t\sg_t^*)^{-1}\nabla f(t,Y_t,Z_t,\wt Z_t)\in H^2(R^n,\sg),
$$
then $Y_t\ge \wt Y_t$ a.s. for all $t\in [0,T]$.

{\it Proof.} Taking the difference of the equations (\ref{eq}),
(\ref{in}) with generators $(f,g,\xi)$ and $(\wt f,\wt g,\wt \xi)$
respectively, we have
$$
Y_t-\wt Y_t= Y_0 -\wt Y_0
$$
$$
-\int_0^t[f(s,Y_s,Z_s)-f(s,\wt Y_s,\wt Z_s)]dK_s
$$
$$
-\int_0^t[f(s,\wt Y_s,\wt Z_s)-\td f(s,\wt Y_s,\wt Z_s)]dK_s-
\int_0^t[g(s)-\td g(s)]d\la N\ra_s
$$
\begin{equation}
-\int_0^t\wt g(s)d(\la N\ra_s-\la \wt N\ra_s)+\int_0^t(Z_s-\wt
Z_s)dM_s+N_t-\wt N_t.
\end{equation}
Let us define the measure $Q$ by $dQ={\cal E}_T(\Lambda)dP$,
where
$$
\Lambda_t=\int_0^t\nabla f(s)^*(\sg_s\sg_s^*)^{-1}dM_s
+\int_0^t\wt g(s)d(N_s+\wt N_s).
$$
By Lemma 3  $Z,\wt Z\in H^2$ and $N$, $\wt N$ are BMO- martingales.
Therefore
Condition $L1),L2)$ and (\ref{2.4}) imply that $\Lambda\in BMO$ and
hence $Q$ is a
probability measure equivalent to $P$.

Denote by $\bar \Lambda$ the martingale part of $\dl Y=Y-\wt Y$, i.e.,
$$
\bar \Lambda=(Z-\wt Z)\cdot M+ N-\wt N.
$$
Therefore, by Girsanov's Theorem and by (\ref{38}) the process
$$
\dl Y_t+\int_0^t(\partial_yf(s)\dl Y_s+\nabla f(s)^*\dl Z_s)dK_s
$$
$$
+\int_0^t\dl f(s,\wt Y_s,\wt Z_s)dK_s+\int_0^t\dl g(s)d\la N\ra_s
$$
$$
=\dl Y_t+\int_0^t(\partial_yf(s)\dl Y_s+\dl f(s,\wt Y_s,\wt
Z_s))dK_s
$$
$$
+\int_0^t\nabla f(s)^*(\sg_s\sg_s^*)^{-1}d\la M\ra_s\dl
Z_s+\int_0^t\dl g(s)d\la N\ra_s
$$
$$=-\int_0^t\wt g(s)d(\la
N\ra_s-\la \wt N\ra_s)+\int_0^t(Z_s-\wt Z_s)dM_s+N_t-\wt N_t $$
$$
=\bar \Lambda_t- \la\Lambda,\bar \Lambda\ra_t,
$$
is a local martingale under $Q$. Moreover, since by Lemma 3 $\bar
N\in BMO$, Proposition 11 of \cite{DdM} implies that
$$
\bar \Lambda_t-\la\Lambda,\bar \Lambda\ra_t\in BMO(Q).
$$
Thus, using the martingale property and the boundary conditions
$Y_T=\xi, \wt Y_T=\wt \xi$ we have
$$
Y_t-\wt Y_t=
$$
$$
=E^Q\big(e^{\int_t^T\partial_yf_sdK_s}(\xi-\td \xi)
$$
$$
+\int_t^T e^{\int_t^s\partial_yf_udK_u}(f(s,\wt Y_s,\wt Z_s)-\td
f(s,\wt Y_s,\wt Z_s)) dK_s|F_t\big)
$$
\begin{equation*}
+E^Q\big(\int_t^Te^{\int_t^s\partial_yf_udK_u}(g(s)-\wt g(s))d\la
N\ra_s|F_t\big),
\end{equation*}
which implies that $Y_t\ge \wt Y_t$ a.s. for all $t\in[0,T]$.\qed

{\bf Corollary.} Let condition A) be satisfied. Then if the
solution of (\ref{eq}),(\ref{in}) exists it is unique.

The proof of {\bf Theorem 1} follows now from the last corollary
and Proposition 2.

{\bf Remark.}  Condition L1),L2) is satisfied if there is constant
$C>0$ such that
$$
|f(t,y,z)-f(t,\td y,\td z)|\le C|y-\td y|+C|z-\td z|(|z|+|\td z|)
$$
and $tr(\sg_t\sg_t^*)^{-1}\le C$ $\text{for all}\;\;\; y, \td y\in
R, \; z, \td z\in R^n\;t\in [0,T].$ Conditions L1),L2) are also
fulfilled if $f(t,y,z)$ satisfies the global Lipschitz condition
and $M\in BMO$.

\

{\bf Acknowledgements}

\

I am thankful to M. Mania for useful discussions and remarks.

I also would  like to thank the referee for valuable remarks and
suggestions which lead to many improvements in the original
version of the paper .

\

\bibliography{alpha}

\end{document}